\def\versiondate{4 April 2020}
\input math10.macros

\let\nobibtex = t
\let\noarrow = t
\input eplain
\beginpackages
\usepackage{url}
   \usepackage{color}  
\endpackages

\input Ref.macros

\checkdefinedreferencetrue
\continuousfigurenumberingtrue
\theoremcountingtrue
\sectionnumberstrue
\forwardreferencetrue
\citationgenerationtrue
\nobracketcittrue
\hyperstrue
\initialeqmacro


\input\jobname.key
\bibsty{../../texstuff/myapalike}

\def\Xsp{X}  

\def\ip#1,#2{\langle #1, #2 \rangle}
\def\Bigip#1,#2{\Big\langle #1, #2 \Big\rangle}
\def\bigip#1,#2{\big\langle #1, #2 \big\rangle}

\def\HH{{\Bbb H}}

\def\as{\bar}
\def\aswide#1{\overline{\smash{#1}\vphantom{\vrule height 5pt}}}
\def\dd{{\rm d}}
\def\bd{\partial}

\def\SRB{\ref b.SRB:measure/, hereinafter referred to as 
{\htmllocref{\bibcode{SRB:measure}}{SRB}}%
\def\SRB{\htmllocref{\bibcode{SRB:measure}}{SRB}}}

\ifproofmode \relax \else\head{To appear in {\it Pacific J. Math.}}
{Version of \versiondate}\fi 
\vglue20pt

\title{Strong Negative Type in Spheres}

\author{Russell Lyons}

\abstract{It is known that spheres have negative type, but only subsets
with at most one pair of antipodal points have strict
negative type. These are 
conditions on the (angular) 
distances within any finite subset of points. We show that
subsets with at most one pair of antipodal points have strong
negative type,
a condition on every probability
distribution of points. This
implies that the function of expected distances to points determines
uniquely the probability measure on such a set. It also implies that the
distance covariance test for stochastic independence, introduced by
Sz\'ekely, Rizzo and Bakirov, is consistent against all alternatives in
such sets.
Similarly, it allows tests of goodness of fit, equality of distributions,
and hierarchical clustering with angular distances.
We prove this by showing an analogue of the Cram\'er--Wold theorem.}

\bottomIII{Primary 
51K99,    
51M10,  
45Q05,  
44A12.  
Secondary
62H20,   
62G20,  	
62H15,   	
62H30.   
}
{Cram\'er--Wold, hemispheres,
expected distances, distance covariance, equality of
distributions, goodness of fit, hierarchical clustering.}
{Research partially supported by NSF grant DMS-1612363.}

\bsection{Introduction}{s.intro}

We introduce the topic by borrowing from \ref b.Lyons:hypneg/.

Let $(X, d)$ be a metric space. One says that $(X, d)$ has \dfn{negative
type} if for all $n \ge 1$ and all lists of
$n$ red points $x_i$ and $n$ blue points $x'_i$ in $X$,
the sum $2\sum_{i, j} d(x_i, x'_j)$ of the distances between the $2n^2$
ordered pairs of points of opposite color is at least the sum $\sum_{i, j}
\big(d(x_i, x_j) + d(x'_i, x'_j)\big)$ of the distances between the
$2n^2$ ordered pairs of points of the same color.
It is not obvious that euclidean space has this property, but it is well
known.
By considering repetitions of $x_i$ and taking limits, we arrive at a
superficially more general property:
For all $n \ge 1$, $x_1, \ldots, x_n \in \Xsp$, and $\alpha_1, \ldots,
\alpha_n
\in \R$ with $\sum_{i=1}^n \alpha_i = 0$, we have 
$$
\sum_{i, j \le n} \alpha_i \alpha_j d(x_i, x_j) \le 0
.
\label e.ntdef
$$
We say that $(\Xsp, d)$ has \dfn{strict negative type} if, for every $n$
and all $n$-tuples of distinct points $x_1, \ldots, x_n$,
equality holds in \ref e.ntdef/ only when $\alpha_i = 0$ for all $i$.
Again, euclidean spaces have strict negative type.
A simple example of a metric space of non-strict negative type is $\ell^1$
on a 2-point space, i.e., $\R^2$ with the $\ell^1$-metric.

A (Borel) probability measure $\mu$ on $X$ has \dfn{finite first moment} if
$\int d(o, x) \,\dd\mu(x) < \infty$ for some (hence all) $o \in X$; write
$P_1(X, d)$ for the set of such probability measures.
Suppose that $\mu_1, \mu_2 \in P_1(X, d)$.
By approximating $\mu_i$ by probability measures of finite support,
we obtain a yet more general property, namely,
that when $\Xsp$ has negative type,
$$
\int d(x_1, x_2) \,\dd(\mu_1 - \mu_2)^2(x_1, x_2) \le 0
.
\label e.negtype
$$
We say that $(\Xsp, d)$ has \dfn{strong negative type} if it has negative
type and equality holds in \ref e.negtype/ only when $\mu_1 = \mu_2$. 
See \refbmulti {Lyons:dcoverrata} for an example of a (countable)
metric space of strict but not strong negative type.
The notion of strong negative type was first defined by \ref b.ZKK/.
\ref b.Lyons:dcov/
used it to show that a metric space $X$ has strong negative type
iff the theory of distance covariance holds in $X$ just as in euclidean
spaces, as introduced by \ref b.SRB:measure/.
Similarly, it allows tests of goodness of fit, equality of distributions,
and hierarchical clustering with angular distances: see the review in \ref
b.SZ:energy/. 
\ref b.Lyons:dcov/ noted that
if $(\Xsp, d)$ has negative type, then
$(\Xsp, d^r)$ has strong negative type when $0 < r < 1$.

Define
$$
a_\mu(x) := \int d(x, x') \,\dd\mu(x')
$$
for $x \in X$ and $\mu\in P_1(X, d)$.
\ref b.Lyons:dcov/ remarked that if $(X, d)$ has negative type, then 
the map $\alpha \colon \mu \mapsto a_\mu$ is
injective on $\mu \in P_1(\Xsp)$
iff $\Xsp$ has strong negative type.
(There are also metric spaces not of negative type for which
$\alpha$ is injective.)

A list of metric spaces of negative type appears as Theorem 3.6 of \ref
b.Meckes:pdms/. All euclidean spaces have strong negative type; see \ref
b.Lyons:dcov/ for a discussion of various proofs.

That real and complex
hyperbolic spaces $\HH^n$ have negative type was shown by \ref
b.Gangolli/, Sec.~4, and was made explicit by \ref b.FarHarz/, Corollary 7.4; 
that they have strict negative type was shown by \ref b.HKM/.
\ref b.Lyons:hypneg/ showed that real hyperbolic spaces have strong
negative type.
The remaining constant-curvature, simply connected spaces are spheres.

Let $S^n$ denote the unit-radius sphere centered at the origin of
$\R^{n+1}$.
Although spheres have negative type (in their intrinsic metric), not even
circles have strict negative type. For example, in $S^1$, take two red
points $\bigl\{(1, 0), (-1, 0)\bigr\}$ and two blue points $\bigl\{(0, 1),
(0, -1)\bigr\}$.
Nevertheless, antipodal symmetry is the only obstruction to strict negative
type: 
the main result, Theorem
9.1, of \ref b.HLMT/ is that
a subset of a sphere has strict negative
type iff that subset contains at most one pair of antipodal points.
We strengthen this to strong negative type:

\procl t.strong-gen
If $B \subset S^n$ contains at most one pair of antipodal points, then $B$
has strong negative type.
\endprocl

We begin by proving a special case:

\procl t.strong
If $H \subset S^n$ is an open hemisphere, then $H$
has strong negative type.
\endprocl

We may parametrize open hemispheres as 
$$
H_t := \{x \in S^n \st t \cdot x > 0\}
$$
for $t \in S^n$.
A crucial ingredient in the proof of \ref t.strong/ is an analogue of the
Cram\'er--Wold theorem:

\procl t.CW
Let $H$ be an open hemisphere in an $n$-dimensional sphere, $S^n$. 
For a finite signed measure $\mu$ on $H$ and $t \in S^n$, define $b_\mu(t)
:= \mu(H \cap H_t)$. The map $\mu \mapsto b_\mu$ is injective. Moreover,
if $D$ is a dense subset of $S^n$, then
$\mu \mapsto b_\mu \restrict D$ is injective.
\endprocl

Let $R\colon x \mapsto -x$ be the reflection in the origin.
If $K$ is a Borel subset of $S^n$ such that $K$ and its image under $R$
partition $S^n$ and such that the interior of $K$ is a hemisphere, then
call $K$ a \dfn{partitioning hemisphere}.
Given a probability measure $\mu$ on $S^n$, let $\mu_R$ denote the maximal
measure that is invariant under $R$ and such that $\mu_R \le \mu$.
Note that if $\mu$ is a probability measure on $S^n$ that is invariant
under $R$, then $\mu(K) = 1/2$ for every partitioning hemisphere, $K$.
Therefore, for every probability measure $\mu$ on $S^n$ with $\mu_R \ne 0$,
there is a probability measure $\nu \ne \mu$ such that $\mu(K) = \nu(K)$
for every partitioning hemisphere, $K$.
Moreover, if $\mu_R \ne 0$ but
$\mu_R(A) = 0$ for every $(n-1)$-dimensional great sphere $A$
in $S^n$, then there is another probability measure $\nu \ne \mu$ such that
$\mu(H) = \nu(H)$ for every open hemisphere, $H$.
We extend \ref t.CW/ to show the converse, which we use to prove
\ref t.strong-gen/:

\procl t.contra
Let $\mu$ be a probability measure on $S^n$ such that $\mu_R = 0$. 
If\/ $\nu$ is a probability on $S^n$ such that $\mu(K) = \nu(K)$ for every
partitioning hemisphere, $K$, then $\nu = \mu$.
Similarly,
if\/ $\nu$ is a probability on $S^n$ such that $\mu(H_t) = \nu(H_t)$ for every
$t$ belonging to a dense subset $D$ of $S^n$, then $\nu = \mu$.
\endprocl

The last assertion of \ref t.contra/ is essentially 
known: see, e.g., Lemmas 2.3 and 2.4 of \ref b.Rubin/.

We also have the following fact:

\procl p.symm
If $\mu$ is a probability measure on $S^n$ 
such that $\mu(K) = 1/2$ for every
partitioning hemisphere, $K$, then $\mu$ is $R$-invariant.
Similarly, 
if $\mu$ is a probability measure on $S^n$ 
such that $\mu(H_t) = \mu(H_{-t})$ for a dense set of
$t$, then $\mu$ is $R$-invariant.
\endprocl

The last assertion is again essentially
known: see Korollar 3.2 of \ref b.Schneider/.
The work of \ref b.Rubin/ and \ref b.Schneider/, as well as other authors
who study related questions, uses spherical harmonics. This is a 
powerful tool that leads to more general results, although those extensions
do not seem relevant to negative type. We give elementary
proofs that rely only on the Cram\'er--Wold theorem for euclidean spaces:

\procl t.euclCW
If $\mu$ is a complex Borel measure on $\R^n$ such that $\mu(H) = 0$ for
every open halfspace in $\R^n$, then $\mu = 0$.
\endprocl

\proof
For $t \in S^{n-1}$, define $\mu_t$ on $\R$ by 
$$
\mu_t(-\infty, a)
:=
\mu\{x \in \R^n \st t \cdot x < a\} \quad (a \in \R).
$$
Then $\mu_t = 0$, whence its Fourier transform $\widehat {\mu_t}$ satisfies
$\widehat {\mu_t}(b) = 0$ for all $b \in \R$. Because $\widehat {\mu_t}(b)
= \widehat \mu(b t)$, it follows that the Fourier transform of $\mu$ also
vanishes, whence so does $\mu$. 
\Qed

Even this theorem can be proved without
Fourier analysis---see \ref b.Walther:withAdd/ or \ref b.LZ/.

\bsection{Proofs}{s.proof}

\proofof t.CW
By the bounded convergence theorem,
$b_\mu \restrict D$ determines $b_\mu(t)$ for all $t$ such that
$\mu\bigl(\partial (H \cap H_t)\bigr) = 0$, and therefore $b_\mu \restrict
D$ determines all of $b_\mu$ by continuity from below: for every $t \in S^n$,
there are $s_k \in S^n$ such that
$\mu\bigl(\partial (H \cap H_{s_k})\bigr) = 0$ and $H \cap H_{s_k}$ increase to
$H \cap H_t$.

We may take $H$ to be the upper open hemisphere, $\bigl\{(t_1, t_2, \ldots,
t_{n+1}) \in S^n \st t_{n+1} > 0\bigr\}$.
Define $\phi\colon H \to \R^{n+1}$ by 
$$
\phi(t_1, \dots, t_{n+1}) := (t_1/t_{n+1}, \dots, t_{n}/t_{n+1}, 1).
$$
Then $\phi$ is a homeomorphism from $H$ to the affine hyperplane $H' :=
\bigl\{(t_1, t_2, \ldots, t_{n+1}) \in \R^{n+1} \st t_{n+1} > 0\bigr\}$,
namely, $\phi(t)$ is the intersection of $H'$ with the line through the
origin and $t$.
Furthermore, $\phi$ maps $H \cap H_t$ to an open halfspace in $H'$ and
every open halfspace in $H'$ is the image under $\phi$ of some $H \cap
H_t$.
Therefore, $b_\mu$ determines the measures of all open
halfspaces with respect to the pushforward $\phi_* \mu$ on $H'$.
The classical theorem of Cram\'er and Wold applied to $H'$ 
shows that this determines $\phi_* \mu$, which in turn determines $\mu$.
\Qed

\proofof t.strong
Write $\sigma$ for the volume measure on $S^n$ normalized to have mass
$\pi$.
Then for all $x_1, x_2 \in S^n$, we have
$$
d(x_1, x_2) 
=
\int \big|\I {H_t}(x_1) - \I {H_t}(x_2)\big|^2
\, \dd\sigma(t)
.
$$
This well-known fact is easy to see: By rotation-invariance of $\sigma$,
the right-hand side depends only on $d(x_1, x_2)$. By considering three
points on a great circle, we find that the dependence is linear. Finally,
by taking antipodal points, we verify that the constant of linearity is 1.

Therefore, if $\mu_1$ and $\mu_2$ are probabilities on $S^n$, we may write
$$
\int d(x_1, x_2) \,\dd(\mu_1 - \mu_2)^2(x_1, x_2) 
=
\int \!\! \int \big|\I {H_t}(x_1) - \I {H_t}(x_2)\big|^2
\,\dd(\mu_1 - \mu_2)^2(x_1, x_2) \, \dd\sigma(t)
.
$$
Expanding the square in the integrand and using the facts that 
$$
\int \I {H_t}(x) \,\dd\nu^2(x, y) = \nu(H_t) \nu(S^n)
$$
and
$$
\int \I {H_t}(x) \I {H_t}(y) \,\dd\nu^2(x, y) = \nu({H_t})^2
$$
for any finite signed measure, $\nu$,
we obtain that
$$
\int d(x_1, x_2) \,\dd(\mu_1 - \mu_2)^2(x_1, x_2) 
=
-2 \int \big(\mu_1({H_t}) - \mu_2({H_t})\big)^2\, \dd\sigma(t)
.
$$
It is evident from this that $(S^n, d)$ has negative type.
In order to prove $(H, d)$ has strong negative type, it suffices to show
that if $\mu_1$ and $\mu_2$ are concentrated on $H$ and satisfy
$\mu_1({H_t}) = \mu_2({H_t})$ for $\sigma$-a.e.\ $t$, then $\mu_1 =
\mu_2$.
But this is immediate from \ref t.CW/.
\Qed

Given any signed measure $\theta$, define the antisymmetric measure
$\as\theta := \theta - R_* \theta$, where $R_* \theta$ is the pushforward of
$\theta$ by $R$.
For positive $\theta$ with $\theta_R = 0$, we have $\theta = \as\theta^+$,
the positive part of $\as\theta$.
For positive $\theta$ without assuming that $\theta_R = 0$, we have
$$
2\,\theta(S^n) \ge |\as\theta|(S^n),
\hbox{ with equality iff } \theta_R = 0.
\label e.compare
$$

\procl l.contra
Let $\mu$ and $\nu$ be probability measures on $S^n$. 
If\/ $\mu(K) = \nu(K)$ for every partitioning hemisphere, $K$, then $\as\mu
= \as\nu$.
Similarly, if\/ $\mu(H_t) = \nu(H_t)$ for every $t \in D$,
where $D$ is a dense subset of $S^n$,
then $\as\mu = \as\nu$.
\endprocl

\proof
We claim that there is an $(n-1)$-dimensional great sphere $A$ in $S^n$
with $\mu(A) = \nu(A) = 0$. To see this, we build $A$ inductively by
dimension. First, because only countably many points have positive mass,
there is a pair $A_0$ of antipodal points with $\mu(A_0) = \nu(A_0) = 0$.
Second, all uncountably many
$1$-dimensional great spheres in $S^n$ that contain $A_0$ have
pairwise intersections exactly $A_0$, whence there is a $1$-dimensional
great sphere $A_1 \supset A_0$ with $\mu(A_1) = \nu(A_1) = 0$. We may
continue this procedure recursively, finding a $k$-dimensional great sphere
$A_k \supset A_{k-1}$ for $1 \le k \le n-1$ with $\mu(A_k) = \nu(A_k) = 0$.
Finally, take $A := A_{n-1}$.

Let $H$ be one of the two open hemispheres comprising $S^n \setminus A$.
Note that $\mu(H) = \nu(H)$ under either assumption (in the second case, we
use a continuity argument like that at the start of the proof of \ref
t.CW/).

Let $K$ be a partitioning hemisphere.
Because $\mu(A) = 0$ and $\mu(S^n) = 1$,
we have
$$
\as\mu(H \cap K) 
=
\mu(H \cap K)
+
\mu(H \cap RK) - \mu(H \cap RK) - \mu(RH \cap RK)
=
\mu(H) - \mu(RK)
=
\mu(H) + \mu(K) - 1.
$$
A similar equation holds for $\nu$.
Hence, the assumption that $\mu(K) = \nu(K)$ for every partitioning
hemisphere, $K$, yields
$$
\as\mu(H \cap K) 
=
\as\nu(H \cap K) 
$$
for every such $K$.

Now every set $H \cap H_t$ is of the form $H
\cap K$ for some partitioning hemisphere, $K$. It follows that 
$$
\as\mu(H \cap H_t) 
=
\as\nu(H \cap H_t) 
$$
for every $t$, whence by \ref t.CW/, it follows that $\as\mu = \as\nu$.

We now prove the second assertion of the lemma.
Note that in the preceding proof, we did not use the full strength of the
assumption that $\mu(K) = \nu(K)$ for every partitioning
hemisphere, $K$, but only that 
$\mu(K_t) = \nu(K_t)$ for partitioning
hemispheres $K_t$ satisfying $K_t \cap H = H_t \cap H$ for $t \in D$; we
may also require that $K_t \cap RH = \overline{H_t} \cap RH$.
Let $u$ be such that $H_u = H$, and let $s_k$ be on the geodesic segment
from $u$ to $t$ with $s_k \notin \{u, t\}$, $\>s_k \to t$, and $\mu(\bd
H_{s_k}) = \nu(\bd H_{s_k}) = 0$.
(Such $s_k$ exist because $\mu(A) = \nu(A) = 0$.)
Let $t_{k, j} \in D$ converge to $s_k$ as $j \to\infty$.
By the bounded convergence theorem, $\lim_{j \to\infty} \mu(H_{t_{k, j}}) =
\mu(H_{s_k})$ and similarly for $\nu$, whence $\mu(H_{s_k}) =
\nu(H_{s_k})$.
In addition, we have
$\lim_{k \to\infty} \mu(H_{s_k}) = \mu(K_t)$ and similarly for $\nu$.
Hence, 
$\mu(K_t) = \nu(K_t)$ for every $t \in D$,
whence $\as\mu = \as\nu$.
\Qed

\proofof t.contra
By \ref l.contra/, either assumption implies that $\as\mu = \as\nu$.
We may conclude from \ref e.compare/ that $2 = 2\, \nu(S^n) \ge
|\as\nu|(S^n) = |\as\mu|(S^n) = 2\, \mu(S^n) = 2$, whence again
from \ref e.compare/, that
$\nu_R = 0$. Since also $\mu_R = 0$,
we obtain the desired conclusion, $\mu = \as\mu^+ = \as\nu^+ = \nu$.
\Qed


\proofof t.strong-gen
If $B$ contains no antipodal points, then every $\mu$ concentrated on $B$ has
$\mu_R = 0$, whence the
proof that $B$ has strong negative type is exactly as for \ref t.strong/,
using \ref t.contra/ in place of \ref t.CW/.

If $B$ contains one antipodal pair, $\{x, Rx\}$, then it still suffices to
show that for probabilities $\mu$ and $\nu$ concentrated on $B$, the
assumption $\mu(H_t) = \nu(H_t)$ for a dense set of $t$ implies $\mu =
\nu$.
By \ref l.contra/, such an assumption
yields $\as\mu = \as\nu$.
Because $\as\mu = \aswide{\mu - \mu_R}$ and $\mu - \mu_R$ is a positive
measure with $(\mu - \mu_R)_R = 0$, and similarly for $\nu$, we obtain
$\aswide{\mu - \mu_R} = \aswide{\nu - \nu_R}$ and
$\mu - \mu_R = \aswide{\mu - \mu_R}^{\,+} = \aswide{\nu
- \nu_R}^{\,+} = \nu - \nu_R$. Therefore, $\mu_R(H_t) =
\nu_R(H_t)$ for a dense set of $t$.
Because $\mu_R$ and $\nu_R$ are supported by
$\{x, Rx\}$, it follows that $\mu_R =
\nu_R$, and so $\mu = \nu$, as desired.
\Qed

\proofof p.symm
For both assertions, we may apply \ref l.contra/ to the pair of
measures $\mu$ and $R_* \mu$, getting
$\as\mu = - \as\mu$, whence $\as\mu = 0$.
Thus, $\mu = R_* \mu$, as desired.
%
\Qed

\medbreak
\noindent {\bf Acknowledgement.}\enspace 
I thank Marcos Matabuena for asking me about strong negative type for
the angular metric on compositional data, i.e., on the probability simplex.

\def\noop#1{\relax}
\input \jobname.bbl

\filbreak
\begingroup
\eightpoint\sc
\parindent=0pt\baselineskip=10pt

Department of Mathematics,
831 E. 3rd St.,
Indiana University,
Bloomington, IN 47405-7106
\emailwww{rdlyons@indiana.edu}
{http://pages.iu.edu/\string~rdlyons/}

\endgroup

\bye